\documentclass{bmcart}

\usepackage{amsthm}
\usepackage{amsfonts}
\usepackage{amssymb}
\usepackage{amsmath}
\usepackage{graphicx}
\usepackage{tikz}
\usetikzlibrary{positioning}
\usepackage{array}
\usepackage{enumerate}
\usepackage{float}
\usepackage{array}
\newcolumntype{P}[1]{>{\centering\arraybackslash}p{#1}}
\usepackage{subcaption}
\usepackage{float}% just for the example
\usepackage[pagewise]{lineno}
\usepackage{wrapfig}
\usepackage{colortbl}

\usepackage{wrapfig}
\usepackage{booktabs}
\usepackage{subcaption}
\usepackage{lipsum}
\usepackage{tikz}
\def\cmark{\tikz\fill[scale=0.4](0,.35) -- (.25,0) -- (1,.7) -- (.25,.15) -- cycle;} 
\usepackage{hyperref}

\usepackage{fourier} 
\usepackage{array}
\usepackage{makecell}

\newtheorem{definition}{Definition}
\newtheorem{example}{Example}

\DeclareFontFamily{U}{cmex}{}{}

\DeclareFontShape{U}{cmex}{m}{n}{
  <-6> cmex5
  <6-7> cmex6
  <7-8> cmex7
  <8-9> cmex8
  <9-10> cmex9
  <10-12> cmex10
  <12-> cmex12}{}{}

\DeclareSymbolFont{Xcmex}{U}{cmex}{m}{n}

\DeclareMathSymbol{\Xdsum}{\mathop}{Xcmex}{88}
\DeclareMathSymbol{\Xtsum}{\mathop}{Xcmex}{80}
\DeclareMathOperator*{\Xsum}{\mathchoice{\Xdsum}{\Xtsum}{\Xtsum}{\Xtsum}}

\startlocaldefs
\endlocaldefs

\begin{document}

%%% Start of article front matter
\begin{frontmatter}

\begin{fmbox}
\dochead{Research}


\title{Persistence Homology of Networks: Methods and Applications}

%%%%%%%%%%%%%%%%%%%%%%%%%%%%%%%%%%%%%%%%%%%%%%
%%                                          %%
%% Enter the authors here                   %%
%%                                          %%
%% Specify information, if available,       %%
%% in the form:                             %%
%%   <key>={<id1>,<id2>}                    %%
%%   <key>=                                 %%
%% Comment or delete the keys which are     %%
%% not used. Repeat \author command as much %%
%% as required.                             %%
%%                                          %%
%%%%%%%%%%%%%%%%%%%%%%%%%%%%%%%%%%%%%%%%%%%%%%

\author[
   addressref={aff1},                   % id's of addresses, e.g. {aff1,aff2}
   corref={aff1},                       % id of corresponding address, if any
                        % id's of article notes, if any
   email={maktas@uco.edu}   % email address
]{\inits{MA}\fnm{Mehmet E} \snm{Aktas}}
\author[
   addressref={aff2},
   email={john.RS.Smith@cambridge.co.uk}
]{\inits{EA}\fnm{Esra} \snm{Akbas}}
\author[
   addressref={aff1},
   email={aelfatmaoui@uco.edu}
]{\inits{AF}\fnm{Ahmed} \snm{El Fatmaoui}}

%%%%%%%%%%%%%%%%%%%%%%%%%%%%%%%%%%%%%%%%%%%%%%
%%                                          %%
%% Enter the authors' addresses here        %%
%%                                          %%
%% Repeat \address commands as much as      %%
%% required.                                %%
%%                                          %%
%%%%%%%%%%%%%%%%%%%%%%%%%%%%%%%%%%%%%%%%%%%%%%

\address[id=aff1]{%                           % unique id
  \orgname{Department of Mathematics and Statistics, University of Central Oklahoma}, % university, etc
 % \street{100 N Drive},                     %
 % \postcode{73034}                                % post or zip code
  \city{Edmond, OK},                              % city
  \cny{USA}                                    % country
}
\address[id=aff2]{%
  \orgname{Department of Computer Science, Oklahoma State University},
%  \street{219 MSCS},
 % \postcode{74078}
  \city{Stillwater, OK},
  \cny{USA}
}

%%%%%%%%%%%%%%%%%%%%%%%%%%%%%%%%%%%%%%%%%%%%%%
%%                                          %%
%% Enter short notes here                   %%
%%                                          %%
%% Short notes will be after addresses      %%
%% on first page.                           %%
%%                                          %%
%%%%%%%%%%%%%%%%%%%%%%%%%%%%%%%%%%%%%%%%%%%%%%

\end{fmbox}% comment this for two column layout

%%%%%%%%%%%%%%%%%%%%%%%%%%%%%%%%%%%%%%%%%%%%%%
%%                                          %%
%% The Abstract begins here                 %%
%%                                          %%
%% Please refer to the Instructions for     %%
%% authors on http://www.biomedcentral.com  %%
%% and include the section headings         %%
%% accordingly for your article type.       %%
%%                                          %%
%%%%%%%%%%%%%%%%%%%%%%%%%%%%%%%%%%%%%%%%%%%%%%

\begin{abstractbox}

\begin{abstract} % abstract
%\parttitle{First part title} %if any
Information networks are becoming increasingly popular to capture complex relationships across various disciplines, such as social networks, citation networks, and biological networks. The primary challenge in this domain is measuring similarity or distance between networks based on topology. However, classical graph-theoretic measures are usually local and mainly based on differences between either node or edge measurements or correlations without considering the topology of networks such as the connected components or holes. In recent years, mathematical tools and deep learning based methods have become popular to extract the topological features of networks. Persistent homology (PH) is a mathematical tool in computational topology that measures the topological features of data that persist across multiple scales with applications ranging from biological networks to social networks. 

In this paper, we provide a conceptual review of key advancements in this area of using PH on complex network science. We give a brief mathematical background on PH, review different methods (i.e. filtrations) to define PH on networks and highlight different algorithms and applications where PH is used in solving network mining problems. In doing so, we develop a unified framework to describe these recent approaches and emphasize major conceptual distinctions. We conclude with directions for future work. We focus our review on recent approaches that get significant attention in the mathematics and data mining communities working on network data. We believe our summary of the analysis of PH on networks will provide important insights to researchers in applied network science.
%\parttitle{Second part title} %if any
%Text for this section.
\end{abstract}

%%%%%%%%%%%%%%%%%%%%%%%%%%%%%%%%%%%%%%%%%%%%%%
%%                                          %%
%% The keywords begin here                  %%
%%                                          %%
%% Put each keyword in separate \kwd{}.     %%
%%                                          %%
%%%%%%%%%%%%%%%%%%%%%%%%%%%%%%%%%%%%%%%%%%%%%%


% MSC classifications codes, if any
%\begin{keyword}[class=AMS]
%\kwd[Primary ]{}
%\kwd{}
%\kwd[; secondary ]{}
%\end{keyword}

\end{abstractbox}
%
%\end{fmbox}% uncomment this for twcolumn layout

\end{frontmatter}

%%%%%%%%%%%%%%%%%%%%%%%%%%%%%%%%%%%%%%%%%%%%%%
%%                                          %%
%% The Main Body begins here                %%
%%                                          %%
%% Please refer to the instructions for     %%
%% authors on:                              %%
%% http://www.biomedcentral.com/info/authors%%
%% and include the section headings         %%
%% accordingly for your article type.       %%
%%                                          %%
%% See the Results and Discussion section   %%
%% for details on how to create sub-sections%%
%%                                          %%
%% use \cite{...} to cite references        %%
%%  \cite{koon} and                         %%
%%  \cite{oreg,khar,zvai,xjon,schn,pond}    %%
%%  \nocite{smith,marg,hunn,advi,koha,mouse}%%
%%                                          %%
%%%%%%%%%%%%%%%%%%%%%%%%%%%%%%%%%%%%%%%%%%%%%%

%%%%%%%%%%%%%%%%%%%%%%%%% start of article main body
% <put your article body there>

%%%%%%%%%%%%%%%%
%% Background %%
%%
%\section{Content}
%Text and results for this section, as per the individual journal's instructions for authors. %\cite{koon,oreg,khar,zvai,xjon,schn,pond,smith,marg,hunn,advi,koha,mouse}

\section{Introduction}
Information networks are important tools to model the relationship between complex data. They exist in multiple disciplines such as social networks, biological networks, the World Wide Web and so on. Analysis of such networks includes many applications such as node classification~\cite{bhagat2011node, akbas2019network}, community detection~\cite{akbas2017truss,Akbas:2017:AGC}, and link prediction~\cite{lopes2010collaboration,sharan2007network}. 

The primary challenge in applied network science is measuring similarity or distance between networks without knowing node correspondences. Since comparing the graphs with the graph isomorphism is computationally expensive \cite{babai2016graph}, many statistically oriented graph similarity measures have been proposed in literature \cite{baur2005network,zager2008graph}. While some of these methods embed the graphs into a feature space and then define distances on that space, other methods define kernel functions on graphs to build similarity measures \cite{vishwanathan2010graph}. Moreover, in graph-theoretic approaches, similarity measures are defined based on the difference in graph-theoretic features such as assortativity, betweenness centrality, small-worldness, and network homogeneity. However, such classical graph-theoretic measures are usually local and mainly based on differences between either node or edge measurements, or correlations without considering the network topology. Therefore, they may have information loss over topological structures, such as the connected components or holes in networks. On the other hand, structural holes in networks can give important information about network topology~\cite{xu_sthole}. For instance, node importance can be measured based on structural holes. The unique characteristics of nodes in the location of structural holes can help to separate the structural holes nodes from other nodes. Moreover, the existence and distribution of structural holes in networks can be used as important topological features for network comparison and classification~\cite{xu2018assessing}.

In recent years, mathematical tools and deep learning based methods have become popular to extract the topological features of networks. Persistent homology (PH) is a mathematical tool in computational topology that measures the topological features of data that persist across multiple scales. Its applications range from biological systems \cite{li2017persistent} to computer vision \cite{adcock2013ring}. The basic idea in PH is to replace the data points with a parametrized family of simplicial complexes, which can roughly be considered as a union of points, edges, triangles, tetrahedron and higher-dimensional polytopes, and encode the change of the topological features (such as the number of connected components, holes, voids) of the simplicial complexes across different parameters for data analysis \cite{ghrist2008barcodes}. For an extensive and rigorous introduction to the computation of persistent homology, we refer readers to the survey papers \cite{patania2017topological, otter2017roadmap}.

Nowadays PH is largely applied for the study of complex networks as a feature extractor since persistent homology gives multi-scales summarization of the graph, unlike the traditional metrics that describe the graph in specific angles. In this paper, we provide a conceptual review of key advancements in the area of using PH on complex network science. 

The paper is structured as follows: In Section \ref{sec:prelim}, we define and give the background on networks, simplicial complex, simplicial homology, and persistent homology. In Section \ref{sec:filt}, we list and compare the filtrations defined for networks. In Section \ref{sec:app}, we highlight different algorithms and applications where PH is used in solving network mining problems. Lastly, we conclude the paper with directions for future work in Section \ref{sec:conc}. 

%%%%%%%%%%%%%%%%%%%%%%%%%%%%%%%%%%%%%%%%%%%%%%
%%                                          %%
%% Backmatter begins here                   %%
%%                                          %%
%%%%%%%%%%%%%%%%%%%%%%%%%%%%%%%%%%%%%%%%%%%%%%

\section{Preliminaries}\label{sec:prelim}

While a network can be represented as a graph, it can also be represented as other \textit{topological} objects. Topology is a branch of mathematics that studies the property of the shapes that are invariant under continuous deformation such as stretching, twisting, bending but not tearing or gluing. For example, a donut and a coffee mug are topologically equivalent since one can transform one to the other continuously. Topological invariants, which are properties of the shapes that do not change under continuous transformation, are useful to detect whether given two shapes are topologically equivalent. The number of connected components, the existence of holes or voids are examples of the topological invariants. Algebraic topology is the area in topology that extracts these invariants of an object by simply counting them or associating algebraic structures, such as vector spaces, to them. For example, for a given topological object $X$, \textit{homology} associates vector spaces $H_i(X)$ for $i=\{0,1,2,...\}$ where the dimension of $H_0(X)$ gives the number of connected components, $H_1(X)$ gives the number of holes, $H_2(X)$ gives the number of voids and so on.    

For a finite set of points, e.g. a point cloud data, homology does not give interesting information. The dimension of $H_0(X)$ gives the number of points, and the dimensions of the higher dimensional homology are zero. It is also similar in a network setting. The dimension of $H_0(X)$ gives the number of disconnected subgraphs, $H_1(X)$ gives the number of loops and the dimensions of the higher dimensional homology are zero since a graph does not have 2 and higher dimensional simplices. Hence, instead of just looking the homology of the finite set of points itself, using (1) a \textit{distance function}, e.g. a correlation or a measure of dissimilarity between points, and (2) a \textit{parameter value}, one can add simplices and check how homology changes across different scales. \textit{Persistent homology} then tracks the change in homology as the parameter value increases and detects which topological features ``persist" across different scales. 

In general, it is very difficult to compute homology of an arbitrary topological object. Hence, instead of doing this, we can approximate a topological object with a \textit{simplicial complex} and then compute the homology, that is actually called \textit{simplicial homology}.

In this section, we define how to define persistent homology on a finite set of points and networks. We first give a formal definition of graphs and explain their characteristics. We then define the simplicial complex and how to compute the simplicial homology. Finally, we briefly explain the persistent homology and two special metrics which are very useful for using persistent homology in data analysis applications. 

\subsection{Graphs}\label{sec:graph}
As a formal definition, a graph $G$ is a pair of sets $G = (V, E)$ where $V$ is the set of vertices and  $E$ is the set of edges of the graph. Networks can be represented via graphs where vertices represent the objects and edges represent the relations between objects. There are different types of graphs to represent different relations between vertices. While in an \textit{undirected} graph, edges link two vertices symmetrically, in a \textit{directed} graph (also called \textit{digraph} in literature), edges link two vertices asymmetrically. If there is a score for the relationship between vertices which could represent the strength of interactions, we can represent this type of relations or interactions by a \textit{weighted} network. In a weighted graph, a weight function $W: E \rightarrow \mathbb{R}$ is defined to assign a weight on each edge. Weights could come from the Euclidean space or other spaces.

A graph $G$ with $n$ vertices can be represented by an $n \times n$ adjacency matrix. Entries of the matrix will be $G_{ij}=1$ for an unweighted graph and $G_{ij}=w_{ij}$ for a weighted graph if there is an edge from vertex $i$ to vertex $j$. If there is no edge between vertex $i$ and vertex $j$, it will be $G_{ij}=0$. 

Furthermore, there are two different graph types we study in this paper. Firstly, a graph is called a \textit{metric graph} if each edge is assigned a positive length and if the graph is equipped with a natural metric where the distance between any two points of the graph (not necessarily vertices) is defined to be the minimum length of all paths from one to the other. Secondly, a graph is called a \textit{dynamic (time-varying) graph} if the graph varies over time, i.e. it can have vertex and edge deletions and additions.

\subsection{Simplicial complex}
Informally, a \textit{simplicial complex} is a topological object which is built as a union of points, edges, triangles, tetrahedron, and higher-dimensional polytopes. The building blocks of a simplicial complex are called \textit{simplices} (plural of \textit{simplex}). Simplices are higher dimensional analogs of points, line segments, and triangles, such as a tetrahedron. We start this section with a formal definition of a simplex.

\begin{figure}[h!]
    \centering
    \includegraphics{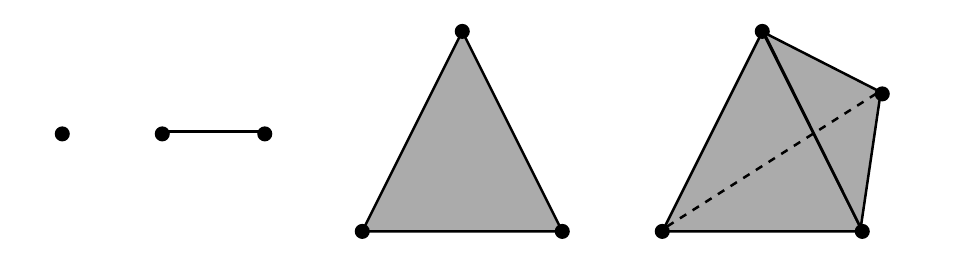}
    \caption{0-,1-,2-, and 3-simplex from left to right}
    \label{fig:simplices}
\end{figure}

\begin{definition}
An \textit{$i$-simplex} $\sigma$ is the convex hull of $i + 1$ affinely independent points, i.e. the set of all convex combinations $\lambda_0 v_0 + \lambda_1 v_1 + ... + \lambda_i v_i$ where $\lambda_0+\lambda_1 + ... + \lambda_i = 1$ and $\lambda_j \leq 1$ for all $j \in \{0,1,...,i\}$.
\end{definition}
A 0-simplex is just a point, a 1-simplex is two points connected with a line segment, a 2-simplex is a filled triangle (see Figure \ref{fig:simplices}). We call \textit{vertex} for 0-simplex, \textit{edge} for 1-simplex, \textit{triangle} for 2-simplex, and \textit{tetrahedron} for 3-simplex.

We can now define a simplicial complex roughly as a union of simplices, but these simplices need to be glued in a certain way. Here is the formal definition.

\begin{definition}\label{def:simpcomp}
A \textit{simplicial complex} is a finite collection of simplices $K$ such that
\begin{enumerate}
    \item Every face of a simplex in $K$ also belongs to $K$.
    \item For any two simplices $\sigma_1$ and $\sigma_2$ in $K$, if $\sigma_1 \cap \sigma_2 \neq \emptyset$, then $\sigma_1 \cap \sigma_2$ is a common face of both $\sigma_1$ and $\sigma_2$.
\end{enumerate}
\end{definition}

\begin{figure}[h!]
    \centering
    \subcaptionbox{ }%
  [.3\linewidth]{\includegraphics[width=0.2\textwidth]{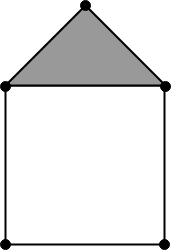}}
\subcaptionbox{ }
  [.3\linewidth]{\includegraphics[width=0.27\textwidth]{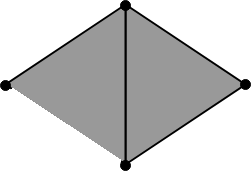}}
  \subcaptionbox{ }
  [.3\linewidth]{\includegraphics[width=0.27\textwidth]{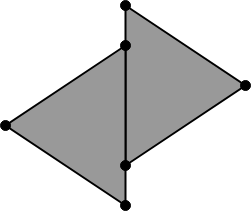}}

    \caption{Finite collection of simplices where (a) is a simplicial complex and, (b) and \\ (c) are not since there is a missing edge in (b) and two triangles meet along an edge\\ which is not an edge
of either triangle in (c).}
    \label{fig:simplex}   
\end{figure}

The first condition says that if a simplex, e.g. a triangle, is in $K$, then its faces, such as its edges and vertices, need to be also in $K$. The second condition says that we can only glue simplices by their common faces. For example, we can glue two triangles by a common vertex or a common edge but cannot glue a vertex of a triangle on one of the edges of the other triangle. Figure \ref{fig:simplex}-a is an example of a simplicial complex whereas Figure \ref{fig:simplex}-b and Figure \ref{fig:simplex}-c are not a simplicial complex since they are violating the first and second condition in Definition \ref{def:simpcomp} respectively.

Before we start to explain how to compute the homology of a simplicial complex, we define the \textit{clique complex} of a graph $G$ which will be a crucial concept to define most of the filtrations in Chapter \ref{sec:filt}. 

\begin{definition}
The clique complex $Cl(G)$ of an undirected graph $G=(V,E)$ is a simplicial complex where vertices of $G$ are its vertices and each $k$-clique, i.e. the complete subgraphs with $k$ vertices, in $G$ corresponds to a $(k-1)$-simplex in $Cl(G)$.   
\end{definition}

\begin{figure*}[h!]
    \centering
    \begin{subfigure}[t]{0.5\textwidth}
        \centering
        \includegraphics[width=.8\textwidth]{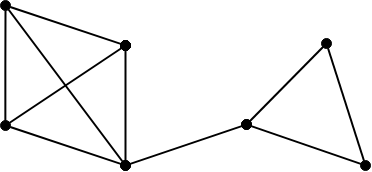}
        \caption{A graph $G$}
    \end{subfigure}%
    ~ 
    \begin{subfigure}[t]{0.5\textwidth}
        \centering
        \includegraphics[width=.8\textwidth]{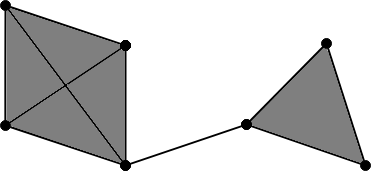}
        \caption{The clique complex $Cl(G)$ of the graph on 
        \\ the left}
    \end{subfigure}
    \caption{An example for constructing the clique complex of a graph}
    \label{fig:clique}
\end{figure*}

For example, in Figure \ref{fig:clique}-a, there is a graph with a 4-clique on the left, 2-clique in the middle and 3-clique on the right. Hence, its clique complex, Figure \ref{fig:clique}-b, has a 3-simplex (tetrahedron), a 1-simplex (edge) and a 2-simplex (triangle). 

\subsection{Simplicial homology}
In a simplicial complex, we can consider the holes as voids bounded by simplices of different dimensions. In dimension 0, they are connected components, in dimension 1, they are loops bounded by edges (1-simplices), in dimension 2, they are holes bounded by triangles (2-simplices) and in general, in dimension $i$, they are the holes bounded by $i$-simplices.

The \textit{simplicial homology} is the way to find the holes in a simplicial complex. To understand what simplicial homology is, we need to define the \textit{chains}, and two special types of chains, namely \textit{cycles} and \textit{boundaries}.

\begin{definition}
Fix a dimension $i$ and assume we use the field of integers. An \textit{$i$-chain} is a formal sum of $i$-simplices of a simplicial complex $K$ with integer coefficients and the sum is taken over possible $i$-simplices. The set of all $i$-chains of $K$ is denoted with $C_i(K)$.
\end{definition}

\begin{wrapfigure}{r}{0.25\textwidth}
  \begin{center}
    \includegraphics[width=0.25\textwidth]{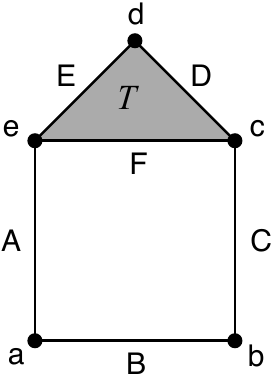}
  \end{center}
  \caption{A simplicial complex with labeled simplices}

\vspace{-.5cm}
  \label{fig:asimp}
\end{wrapfigure}

For example, $c_1=a-3b+4d, c_2=2a+c-2d+3e$ are 0-chains for the simplicial complex in Figure \ref{fig:asimp}. One can add two $i$-chains by simply adding the corresponding integer coefficients, e.g. $c_1+c_2=3a-3b+c+2d+3e$, and multiply by scalars, e.g. $2c_1=2a-6b+8d$. Hence, $C_i(K)$ is actually a vector space over integers (more generally we can over any field such as real numbers). For simplicity, we assume the field is the binary field $\mathbb{Z}/2\mathbb{Z}=\{0,1\}$ from now on.

To map an $i$-simplex to an $(i-1)$-simplex, we define the \textit{boundary} of an $i$-simplex as the sum of its $(i-1)$-dimensional faces. Formally speaking, for an $i$-simplex $\sigma=[v_0,...,v_i]$, its boundary is 
$$
\partial_i\sigma=\Xsum_{j=0}^i [v_0,...,\hat{v_j},...,v_i]
$$
where the hat indicates the $v_j$ is omitted. We can expand this definition to $i$-chains. For an $i$-chain $c=c_i\sigma_i$, $\partial_i(c)=\Xsum c_i\partial_i \sigma i$. For example, in Figure \ref{fig:asimp}, $\partial_1 \text{A}=\text{e}+\text{a}$ and $\partial_2 T = \text{E}+ \text{D} +\text{F}$.

We should also note here that the boundary of a boundary is empty, i.e. $\partial_i\partial_{i+1}=0$. For example, in Figure \ref{fig:asimp}, $\partial_1\partial_2(T)=\partial_1(E+D+F)=(d+e)+(e+c)+(c+d)= 2c+ 2d + 2e = 0$ since $2=0$ in the binary field $\mathbb{Z}/2\mathbb{Z}=\{0,1\}$. 

We can now distinguish two special types of chains using the boundary map that will be useful to define homology. The first one is an \textit{$i$-cycle}, which is defined as an $i$-chain with empty boundary. In other words, an $i$-chain $c$ is an $i$-cycle if and only if $\partial_i(c)=0$, i.e. $c\in \text{ker}(\partial_i)$. For example, the 1-chain $A+B+C+F$ in Figure \ref{fig:asimp} is a 1-cycle since $\partial_1(A+B+C+F)=\partial_1(A)+\partial_1(B)+\partial_1(C)+\partial_1(F)=(e+a)+(a+b)+(b+c)+(c+e)=0$. The set of all such $i$-cycles forms a subspace in $C_i(K)$, which we denote $Z_i(K)$.  

Second special type of an $i$-chain is \textit{$i$-boundary}: an $i$-chain $c$ is an $i$-boundary if there exists an $(i+1)$-chain $d$ such that $c=\partial_{i+1}(d)$, i.e. $c \in \text{im}(\partial_{i+1})$. For example. the one chain $E+D+F$ is a 1-boundary since $E+D+F=\partial_2(T)$. The set of all such $i$-boundaries forms a subspace in $C_i(K)$, which we denote $B_i(K)$. 

After defining these two special subspaces, $i$-cycles $Z_i(K)$ and $i$-boundaries $B_i(K)$ of $C_i(K)$, we now take the quotient space of $B_i(K)$ as a subset of $Z_i(K)$. In this quotient space, there are only the $i$-cycles that do not bound an $(i+1)$-complex left. These are actually the $i$-voids of $K$. We call this quotient space as the \textit{$i$-th homology} of the simplicial complex $K$
$$
\displaystyle H_i(K)=\frac{Z_i(K)}{B_i(K)}=\frac{\text{ker}(\partial_i)}{\text{im}(\partial_{i+1})}.
$$
The dimension of $i$-th homology is called the \textit{$i$-th Betti number} of $K$, $\beta_i(K)$, where $\beta_i(K)= \text{ dim ker }(\partial_i) - \text{ dim im }(\partial_{i+1})$. Basically, the $i$-th Betti number is the number of $i$-dimensional voids in the simplicial complex. For example, $\beta_0$ gives the number of connected components and $\beta_1$ gives the number of loops. In Figure \ref{fig:asimp}, $\beta_0=1, \beta_1=1, \beta_2=0$. 

\subsection{Persistent homology}

\begin{figure}[h!]
    \centering
    \includegraphics[width=.94\textwidth]{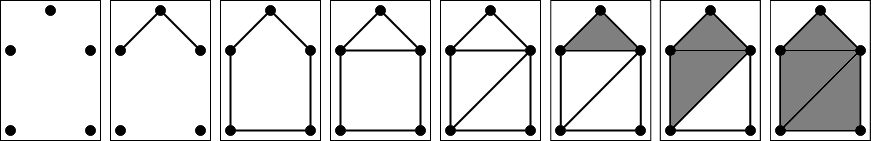}
    \caption{A filtration for $\delta=0,1,2,3,4,5,6,7$ (from left to right)}
    \label{fig:filter}
\end{figure}

For a finite set of points $X$, e.g. a point cloud data, homology does not give interesting information. $\beta_0$ gives the number of connected components, which is just the number of points, and all other Betti numbers are zero since there are no other dimensional holes in the set. Hence, instead of working with the set of points, one can induce a family of simplicial complexes $K_X^{\delta}$ for a range of values of $\delta \in \mathbb{R}$ out of the set $X$ so that the complex at step $m$ is embedded in the complex at $n$ for $m\leq n$, i.e. $K_X^m \subseteq K_X^n$. This nested family of simplicial complexes is called \textit{\textbf{filtration}} (see Figure \ref{fig:filter} for an example). During this construction, some holes may appear and then disappear and the \textit{persistency} of these homological features can be considered as the features of the dataset. In a filtration, one can record the \textit{birth}, the time a hole appears, and \textit{death}, the time a hole disappears, of holes. The essence of the \textit{persistent homology} is to tract the birth and death of these homological features in $K_X^{\delta}$ for different $\delta$ values. The lifespan of each homological feature can be represented as an interval, where the start and end points of the interval correspond to the birth and the death of the homological feature respectively. For a given dataset and a filtration, one can record all these intervals by a \textit{persistence barcode} (PB) as a multiset of intervals bounded below \cite{carlsson2005persistence}. Equivalently, a persistence barcode can be represented via \textit{persistence diagram} (PD) that consists of the birth and death times of the features as a point (birth, death) in the extended real plane $\bar{\mathbb{R}}^2$ \cite{edelsbrunner2000topological}. The longer bars in PBs and the points far away to diagonal in PDs are considered as the real feature of the dataset. 

\begin{figure*}[h!]
    \centering
    \begin{subfigure}[t]{0.5\textwidth}
        \centering
        \includegraphics[width=.8\textwidth]{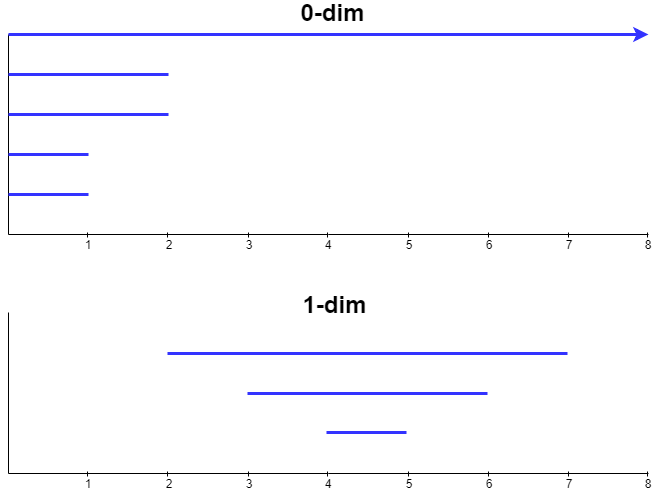}
        \caption{}
    \end{subfigure}%
    ~ 
    \begin{subfigure}[t]{0.5\textwidth}
        \centering
        \includegraphics[width=.7\textwidth]{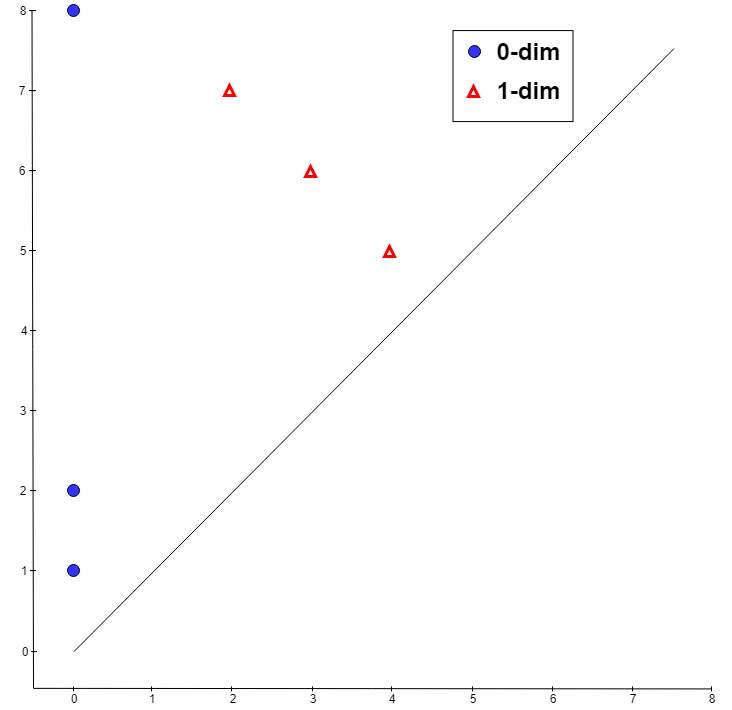}
        \caption{}
    \end{subfigure}
    \caption{Results of the filtration in Figure \ref{fig:filter}, (a) Persistence barcodes for 0- and 1-dim \\(b) Persistence diagrams for 0- and 1- dim }
    \label{fig:pd}
\end{figure*}

\begin{example}
Figure \ref{fig:pd} has the 0- and 1-dimensional persistence barcodes (Figure \ref{fig:pd}-(a)) and the 0- and 1-dimensional persistence diagrams (Figure \ref{fig:pd}-(b)) of the filtration in Figure \ref{fig:filter}. 

We first investigate the 0-dimensional PB and PD. As we see in the filtration, when $\delta=0$, there are five disconnected vertices, which means there are five connected components in the simplicial complex. That is why five bars are \textit{born} at the beginning of the 0-dimensional PB. When $\delta=1$, two edges are added that decreases the total number of connected components to three, hence two bars \textit{die} at $\delta=1$. When $\delta=2$, three more edges are added and this makes the simplicial complex only one connected component, thus two more bars \textit{die} at $\delta=2$. After this point, the number of connected component does not change so the top bar lives forever (arrowhead at the right of that bar implies this fact). Following the same reasoning, 0-dimensional PD has the point $(0,1)$ two times, that corresponds to the two bars spanning from 0 to 1 in the 0-dimensional PB, the point $(0,2)$ again two times, that corresponds to the two bars spanning from 0 to 2 in the 0-dimensional PB, and the point $(0,\infty)$, that corresponds to the top bar that lives forever in the 0-dimensional PB.

For the 1-dimensional PB and PD, since the first 1-dimensional hole (loop) appears for $\delta=2$, there is a bar born at this value in the 1-dimensional barcode. When $\delta=3$, this loop splits into two loops, hence the number of the loops increases to two, and as a result, a new bar is born at $\delta=3$. When $\delta=4$, one of the two loops also splits into two loops, so there is another bar born at $\delta=4$. When $\delta=5$, the top triangle is filled in (a 2-simplex is added), so the number of the loops decreases by one and this results into a death of the bar born at $\delta=4$. Similarly, other two bars die at $\delta=6$ and $\delta=7$. In the 1-dimensional PD, there are the points (2,7), (3,6) and (4,5) that correspond to the three bars in the 1-dimensional PB.
\end{example}

\subsection{Two metrics for persistence diagrams}

One may want to employ persistence diagrams to compare the corresponding datasets. For example, in the network matching problem, we can create a persistence diagram for each network and compare the persistence diagrams to obtain the network similarity. For such a comparison, we need to measure the distance between persistent diagrams using \textit{stable} metrics. A metric is \textit{stable} if a small perturbation of a dataset creates only a small change in the persistence diagram up to that metric. There are two metrics, which can be stable depending on how simplices are defined, that have been commonly used to measure the distance between diagrams: the bottleneck distance and the Wasserstein distance. We first define the bottleneck distance. 

\begin{definition}
Let $P$ and $Q$ be two persistence diagrams. The \textit{bottleneck distance} between $P$ and $Q$ is defined as 
$$
d_B(P,Q)=\inf_{\gamma}\sup_{x\in P}||x-\gamma(x)||_{\infty},
$$
where $\gamma$ ranges over all matchings from $P$ to $Q$ and $||p-q||_{\infty}=\max(|p_1-q_1|,|p_2-q_2|)$ for $p=(p_1,p_2), q=(q_1,q_2) \in \bar{\mathbb{R}^2}$ with $|\infty - \infty|=0$. 
\end{definition}

In other words, the bottleneck distance measures the distance between two persistence diagrams $P$ and $Q$ by the maximum distance between two points in a matching from $P$ to $Q$. Hence, the bottleneck distance only outputs the distance between the greatest outlier, rather than the distance between all pair of points. 

As an answer to this concern, the Wasserstein distance can be used.

\begin{definition}
Let $P$ and $Q$ be two persistence diagrams. The \textit{$p$-th Wasserstein distance} between $P$ and $Q$ is defined as 
$$
d_{W_p}(P,Q)=\inf_{\gamma}\left(\Xsum_{x\in P}||x-\gamma(x)||_p\right)^{1/p},
$$
where $\gamma$ ranges over all matchings from $P$ to $Q$ and $||p-q||_{p}=(|p_1-q_1|^p+|p_2-q_2|^p)^{1/p}$ for $p=(p_1,p_2), q=(q_1,q_2) \in \bar{\mathbb{R}^2}$ with $|\infty - \infty|=0$. 
\end{definition}

In other words, the Wasserstein distance considers the total distance between the matched pair of points, hence provides an overall quantification for the similarity between persistence diagrams.

\section{Filtrations}\label{sec:filt}
In this section, we review the filtrations defined for networks in the literature. We compare the filtrations according to their properties such as sensitivity to different network types (e.g. directed/undirected, weighted/unweighted). We also provide a comparison table, Table ~\ref{table:conf}, at the end of the section. 

Throughout this section, we use the notations defined for graphs in Section \ref{sec:graph}. 

\hypertarget{vr}{\subsection{Vietoris-Rips filtration (VR)}} Let $G=(V,E)$ be an undirected weighted graph with the weight function $W: V \times V \rightarrow \mathbb{R}$ defined on $E$. For any $\delta \in \mathbb{R}$, the 1-skeleton $G_{\delta}=(V_{\delta},E_{\delta}) \subset G$ is defined as the subgraph of $G$ where $V_{\delta}= V$ and its edge set $E_{\delta} \in E$ only includes the edges whose weight is less than or equal to $\delta$. Then, for any $\delta \in \mathbb{R}$, we define the Vietoris-Rips complex as the clique complex of the 1-skeleton $G_{\delta}$, $Cl(G_{\delta})$, and the \textit{Vietoris-Rips filtration} is then defined as $$
\{Cl(G_{\delta}) \hookrightarrow Cl(G_{\delta'}) \}_{0 \leq \delta \leq \delta'}.
$$

In other words, in this filtration, we first start with the vertex set. We then rank the edge weights from the minimum weight, $w_{min}$, to the maximum weight, $w_{max}$, and let the parameter $\delta$ increase from $w_{min}$ to $w_{max}$. At each step, we add the corresponding edges and take the clique complex of the thresholded subgraph $G_{\delta}$. This construction yields the Vietoris-Rips filtration on networks.

For the application purposes, we may prefer to add edges with larger weights before the ones with smaller weights to stress the importance of the weights. In other words, after adding the vertex set, we rank the edge weights from $w_{max}$ to $w_{min}$ and for any $\delta \in \mathbb{R}$, we add edges whose weight is bigger than or equal to $\delta$. This yields a similar yet another filtration. This filtration is called the weight rank clique filtration by \cite{petri2013topological}. However, to be more concise, we prefer to call this as \textit{inverse Vietoris-Rips filtration}.

\hypertarget{dss}{\subsection{Dowker sink and source filtration (DSS)}} Using the idea of the Vietoris-Rips filtration, the authors in \cite{chowdhury2016persistent,chowdhury2018persistent} define the \textit{Dowker $\delta$-sink and $\delta$-source simplicial complex} on directed weighted networks that is sensitive edge directions. For a directed graph $G=(V,E)$ with edge weights $W: V \times V \rightarrow \mathbb{R}$, the Dowker $\delta$-sink simplicial complex associated to $G$ is defined as
$$\mathfrak{D}_{\delta,G}^{si}:=\{ \sigma =[x_0,...,x_n]: \text{ there exists } x'\in V \text{ such that } W(x_i,x')\leq \delta \text{ for each } x_i \}.$$
In other words, there is a \textit{sink} vertex $x'\in V$ such that there are edges from each $x_i \in \sigma$ to $x'$ with weights less than or equal to the threshold $\delta$. Using this simplicial complex, they define the \textit{Dowker sink filtration} as follows
$$
\{\mathfrak{D}_{\delta}^{si} \hookrightarrow \mathfrak{D}^{si}_{\delta'} \}_{0 \leq \delta \leq \delta'}.
$$

They similarly define a dual construction, namely the \textit{Dowker $\delta$-source simplicial complex} associated to a directed weighted network $G$, as follows
$$\mathfrak{D}_{\delta,G}^{so}:=\{ \sigma =[x_0,...,x_n]: \text{ there exists } x'\in V \text{ such that } W(x',x_i)\leq \delta \text{ for each } x_i \}.$$
The only difference here is the edge directions: there is a \textit{source} vertex $x'\in V$ such that there are edges from $x'$ to each $x_i \in \sigma$ with weights less than or equal to the threshold $\delta$. Similarly, they define the \textit{Dowker source filtration} as
$$
\{\mathfrak{D}_{\delta}^{so} \hookrightarrow \mathfrak{D}^{so}_{\delta'} \}_{0 \leq \delta \leq \delta'}.
$$

Dowker sink and source (DSS) filtrations are formed with respect to a central authority $x' \in V$, hence they could be preferred on networks, such as small-world networks, who would desire simplices to be formed with respect to particular hub nodes.

The authors also prove that both filtrations generate the same persistent diagram.

\hypertarget{ccl}{\subsection{Clique complex filtration (CCL)}}  
    
For a graph $G$ with $n \in \mathbb{Z}$ vertices and its clique complex $Cl(G)$, the \textit{clique complex filtration} is defined as 
$$ Cl_0(G) \hookrightarrow Cl_1(G) \hookrightarrow \cdots  \hookrightarrow Cl_n(G) $$ 
such that $Cl_0(G) \subset Cl_1(G) \subset \cdots \subset Cl_n(G)=Cl(G)$ where the $i$-th complex in the filtration is given by $Cl_i(G)=\Xsum_{j=1}^i S_j$ where $S_j$ is the $j$th skeleton of the clique complex, i.e. the set of simplices of dimension less than or equal to $j$ \cite{horak2009persistent}. In other words, in this filtration, we add the vertices at $\delta=0$, add the edges at $\delta=1$, add the triangles at $\delta=2$ and so on.

\hypertarget{vbcl}{\subsection{Vertex-based clique filtration (VBCL)}} This filtration is originally defined in \cite{Rieck2018CliqueCP} for just the 0-dimension, however it can be extended to higher dimensions as well. Let $G=(V,E)$ be a graph with a vertex weight function $\omega: V \rightarrow \mathbb{R}$. For this filtration, we use vertex weights, instead of the edge weights, as threshold values. For any $\delta \in \mathbb{R}$, the 1-skeleton $G_{\delta}=(V_{\delta},E_{\delta}) \subset G$ is defined as the subgraph of $G$ where $V_{\delta}:=\{v \in V | \omega(v)\leq \delta\}$ and the edges $E_{\delta}:=\{e=\{u,v\}\in E | \omega(e):=\text{ max }(\omega(u),\omega(v))\leq \delta \}.$ Then, for any $\delta \in \mathbb{R}$, using the clique complex of the 1-skeleton $G_{\delta}$, $Cl(G_{\delta})$, the filtration is defined as 
$$
\{Cl(G_{\delta}) \hookrightarrow Cl(G_{\delta'}) \}_{0 \leq \delta \leq \delta'}.
$$

Furthermore, we can also define the \textit{inverse vertex-based clique filtration} by just filtering from $\omega_{max}$ to $\omega_{min}$ as we do in the Vietoris-Rips filtration.
  
\hypertarget{wscl}{\subsection{$k$-clique filtration (kCL)}} This filtration is used to detect the evolution of $k$-clique communities for a fixed $k$ in \cite{Rieck2018CliqueCP}. In this filtration, we assume the graph $G=(V,E)$ has a vertex weight function $\omega: V \rightarrow \mathbb{R}$. First, using the vertex weights, we assign a weight function $\omega(\cdot)$ on an arbitrary clique $\sigma$ inductively as 
$$
\omega(\sigma):= \{\smash{\displaystyle\max_{v \subseteq \sigma}} \text{ } \omega(v) | v \text{ is a vertex of the clique } \sigma\},
$$
i.e., the maximum weight of its vertices. Second, for a fixed $k$, we detect all $k$-clique communities in $G$ and create the \textit{$k$-clique connectivity} graph $G^k=(V^k, E^k)$ where there is a vertex for every $k$-clique of $G$ and its edges are defined by 
$$
E^k:=\{ (\sigma,\sigma') \in V^k \times V^k | \sigma \text{ and } \sigma' \text{ are adjacent}\}, 
$$
i.e., $\sigma$ and $\sigma'$ intersect in a $(k-1)$-clique, in other words, they share $k-1$ vertices in common. We then extend the weight function $\omega(\cdot)$ to the edges of $G_k$ by setting
$$
\omega(\sigma, \sigma'):= \max(\omega(\sigma),\omega(\sigma')). 
$$

Next, in a similar way, for any $\delta \in \mathbb{R}$, the 1-skeleton $G^k_{\delta}=(V^k_{\delta},E^k_{\delta}) \subset G^k$ is defined as the subgraph of $G^k$ where $V^k_{\delta}:=\{v \in V^k | \omega(v)\leq \delta\}$ and the edges $E^k_{\delta}:=\{e=\{u,v\}\in E^k | \omega(e)\leq \delta \}.$ Then, for any $\delta \in \mathbb{R}$, using the clique complex of the 1-skeleton $G^k_{\delta}$, the filtration is defined as 
$$
\{Cl(G^k_{\delta}) \hookrightarrow Cl(G^k_{\delta'}) \}_{0 \leq \delta \leq \delta'}.
$$
This filtration is unique in a sense that it just focuses on the evolution of the $k$-clique communities only in the original graph.

\hypertarget{wscl}{\subsection{Weighted simplex filtration (WS)}} In the previous filtration, we assign weights to arbitrary cliques, i.e. simplices, using the vertex weights. Alternatively, one may use another way to assign weights to simplices and use these weights to create a filtration. For example, Huang et all~\cite{huang2017persistent} assign weights to each simplex in a simplicial complex $K$ based on relationship functions in a given dissimilarity network. For any $\delta \in \mathbb{R}$, they define $K_{\delta} \subseteq K$ to be the collection of simplices appearing before or on $\delta$. Then, this construction yields the filtration
$$
\{K_{\delta}) \hookrightarrow K_{\delta'}) \}_{0 \leq \delta \leq \delta'}.
$$
To be a well-defined filtration, we need to have all faces of each simplex and intersections of any simplices in $K_{\delta}$ also appear before or on $\delta$. They prove that this filtration from a given dissimilarity network is a well-defined filtration, i.e. satisfies both conditions.

\hypertarget{vfb}{\subsection{Vertex function based filtration (VFB)}}
Let $G=(V,E)$ be an undirected graph and $f:V \rightarrow \mathbb{R}$ be a function defined on its vertices. We construct the \textit{sublevel graphs} $G_{\delta}=(V_{\delta},E_{\delta})$ for $\delta \in \mathbb{R}$ where $V_{\delta}=\{ v \in V: f(v)\leq \delta$ and $E_{\delta}=\{ (v_1,v_2)\in E : v_1,v_2 \in V_{\delta}\}$. Hence, increasing $\delta$ from $-\infty$ to $\infty$ provides a nested sequence of increasing subgraphs. The \textit{sublevel vertex function based (VFB)} filtration is given by taking the clique complex of each sublevel graph 
$$
\{Cl(G_{\delta}) \hookrightarrow Cl(G_{\delta'}) \}_{0 \leq \delta \leq \delta'}.
$$

Similarly, we can construct the \textit{superlevel graphs} $G_{\delta}=(V_{\delta},E_{\delta})$ for $\delta \in \mathbb{R}$ where $V_{\delta}=\{ v \in V: f(v)\geq \delta$ and $E_{\delta}=\{ (v_1,v_2)\in E : v_1,v_2 \in V_{\delta}\}$. This time decreasing $\delta$ from $\infty$ to $-\infty$ provides a nested sequence of increasing subgraphs which yields to the \textit{superlevel vertex function based (VFB)} filtration as follows
$$
\{Cl(G_{\delta}) \hookrightarrow Cl(G_{\delta'}) \}_{\delta \geq \delta'}.
$$

\hypertarget{ic}{\subsection{Intrinsic \v{C}ech filtration (IC)}} This filtration is defined only for metric graphs in \cite{gasparovic2017complete}. Let $(G,d_G)$ be a metric graph with geometric realization $|G|$. For any point $x \in |G|$, we define the set $B(x,\delta) := \{y \in |G|: d_G(x,y) < \delta \}$, and we let $U_{\delta}:= \{B(x,\delta):x \in |G|\}$ be an open cover. Since $|G|$ has all the vertices and every point along the edges, it has uncountable points. Hence, $U_{\delta}$ is also an uncountable cover. We let $C_{\delta}$ denote the nerve of $U_{\delta}$ where the nerve of a family of sets $\{Y_i\}_{i \in I}$ is the abstract simplicial complex defined on the vertex set $I$ where a family $\{i_0, i_1,...,i_k\}$ with $i_j \in I$ for all $j\in \{0,...,k\}$ spans a $k$-simplex if and only if $U_{i_0} \cap U_{i_1} \cap \cdots \cap U_{i_k} \neq \emptyset$. The associated intrinsic \v{C}ech filtration is defined as the set of inclusion maps
$$
\{C_{\delta} \hookrightarrow C_{\delta'} \}_{0 \leq \delta \leq \delta'}.
$$

\hypertarget{fmg}{\subsection{Functional metric graph filtration (FMG)}} This filtration is also defined for metric graphs only \cite{dey2015comparing}. Let $G$ be a metric graph and take a fixed point $s \in G$. They consider $f:G \rightarrow \mathbb{R}$ where $f(x)=d_G(x,s)$, i.e. to be the geodesic distance from $x$ to $s$. Let $G_{\delta}:=\{x \in G | f(x)\geq \delta \}$ denote the super-level set of $G$ with respect to $\delta \in \mathbb{R}$. Clearly $G_{\delta} \subseteq G_{\delta'}$ for $\delta \geq \delta'$. Then the filtration is given by 
$$
\{G_{\delta} \hookrightarrow G_{\delta'} \}_{\delta \geq \delta'}.
$$

Similarly, instead of using super-level set of $G$, one can use the sub-level set of $G$ for each $\delta$, $G_{\delta}:=\{x \in G | f(x)\leq \delta \}$. This yields another filtration
$$
\{G_{\delta} \hookrightarrow G_{\delta'} \}_{0\leq \delta \leq \delta'}.
$$
Depending on problems, one may choose either of the filtrations. 

\hypertarget{pow}{\subsection{Power filtration (POW)}} Let $G=(V,E)$ be a graph. A $u-v$ walk is an alternating sequence of vertices and edges beginning with $u$ and ending with $v$ such that every edge joins the vertices immediately preceding and following it. A $u-v$ path is a $u-v$ walk in which no vertex is repeated and the number of edges it contains is its length. The graph distance $d(u,v)$ between $u,v \in V(G)$ is the minimum length of all $u-v$ paths. One can also consider the edge weights while computing the graph distance as well. The $r$th power $G^r,r\geq1$, of $G$ is the graph with vertex set $V(G^r)=V(G)$ and for which $\{u,v\}\in E(G^r)$ if, and only if, the distance between $u$ and $v$ in $G$ is at most $r$. The power filtration is the clique complex of the $r$th power $G^r$. In other words, for an appropriate distance range $1\leq r \leq p$ within $G$, the power filtration is given by
$$
Cl(G) \hookrightarrow Cl(G^2) \hookrightarrow \cdots \hookrightarrow Cl (G^p).
$$
where $Cl$ denotes the clique complex.

\hypertarget{tmp}{\subsection{Temporal filtration (TMP)}} This filtration is defined in \cite{Siddharth2017} for dynamic (time-varying) networks. If the network is growing in time $t_0 < \cdots < t_n$, this will yield a sequence of networks $\{G_t, t=t_0,...,t_n\}$ where the network $G_t$ represents the network occurred until time $t$. This network sequence results in the temporal filtration given by the clique complex of each network
$$
Cl(G_{t_0}) \hookrightarrow  Cl(G_{t_1}) \hookrightarrow \cdots \hookrightarrow Cl(G_{t_n}).
$$

\hypertarget{zsf}{\subsection{Zigzag simplicial filtration (ZSF)}} This filtration is also defined for dynamic networks. While \hyperlink{tmp}{TMP} only considers the vertex and edge insertion into a dynamic graph which yields adding simplices to the simplicial complexes, this method also allows vertex and edge deletion from a dynamic graph which yields removing simplices from the simplicial complexes. In a standard filtration on a graph $G$, $G_{\delta} \subseteq G_{\delta'}$ whenever $\delta \leq \delta'$. This filtration generalizes standard filtrations by allowing the simplicial complexes to sometimes become smaller. A zigzag simplicial filtration on a graph $G$ is a filtration with extra two conditions: (1) The set of points of discontinuity of the zigzag simplicial filtration should be locally finite, i.e. each point in the set has a neighborhood that includes only finitely many of the points in the set and (2) for any scale parameter value $\delta \in \mathbb{R}$, it holds that $G_{\delta - \epsilon} \subseteq G_{\delta} \supseteq G_{\delta+\epsilon}$ for all sufficiently small $\epsilon >0$. Then we use the zigzag persistent homology to obtain the persistence barcodes/diagrams \cite{carlsson2010zigzag}. The same basic idea applies in the zigzag persistent homology. For example, for the 0-dimensional zigzag persistence barcodes, we just track the number of connected components in the filtration. 

\hypertarget{pph}{\subsection{Digraph filtration using Persistent Path Homology (PPH)}} 
In \cite{chowdhury2018persistent2}, the authors define a new way to construct homology on networks: persistent path homology (PPH) which is sensitive to the edge directions. We summarize the construction in 4 steps.

\textbf{Step 1:} Let $G=(V, E)$ be directed weighted graph. Given any integer $p \in \mathbb{Z}_+$, an \textit{elementary $p$-path over} $V$ is a sequence $[v_0, ... , v_p]$ of $p+1$ vertices of $V$. For each $p \in \mathbb{Z}_+$, the free vector space consisting of all formal linear combinations of elementary $p$-paths over $V$ with coefficients in $\mathbb{K}$ is denoted $\Lambda_p = \Lambda_p(V) = \Lambda_p(V, \mathbb{K})$. One also defines $\Lambda_{-1} := \mathbb{K}$ and $\Lambda_{-2} := \{0\}$. Next, for any $p \in \mathbb{Z}_+$, one defines the \textit{non-regular boundary map} $\partial_p^{nr}:\Lambda_p \rightarrow \Lambda_{p-1}$ as 
$$
\partial_p^{nr}([v_0,...,v_p]):= \Xsum_{i=0}^p(-1)^i[v_0,...,\hat{v_i}, ..., v_p]
$$
for each elementary $p$-path $[v_0,...,v_p] \in \Lambda_p$. $\partial_{-1}^{nr}: \Lambda_{-1} \rightarrow \Lambda^{-2}$ is the zero map. Observe that $\partial_{p+1}^{nr} \circ \partial_{p}^{nr}=0$ for all $p \geq -1$ so $(\Lambda_p,\partial_p^{nr})_{p\in \mathbb{Z}^+}$ is a chain complex.

\textbf{Step 2:} For each $p\in \mathbb{Z}^+$, an elementary $p$-path $[v_0,...,v_p]$ is called \textit{regular} if $v_i\neq v_{i+1}$ for each $0 \leq i \leq p-1$ and irregular otherwise. Let $\mathcal{R}_p(V,\mathbb{K}):= \mathbb{K}[\{[v_0,...,v_p]: [v_0,...,v_p] \text{ is regular } \}]$ and $\mathcal{I}_p$ is irregular one. We have $\partial_p^{nr}(\mathcal{I}_p) \subseteq \mathcal{I}_{p-1}$ so $\partial_p^{nr}$ is well defined on $\Lambda_p/\mathcal{I}_p$. Since $\mathcal{R}_p \cong \Lambda_p/\mathcal{I}_p$ via a natural isomorphism, one can define $\partial_p: \mathcal{R}_p \rightarrow \mathcal{R}_{p-1}$ as the pullback of $\partial_p^{nr}$ via this isomorphism. $\partial_p$ is called the \textit{regular boundary map} and now we have a chain complex $(\mathcal{R}_p, \partial_p)_{p\in \mathbb{Z}^+}$.

\textbf{Step 3:} For each $p \in \mathbb{Z}^+$, one defines an elementary $p$-path $[v_0,...,v_p]$ on $V$ to be \textit{allowed} if $(v_i,v_{i+1}) \in E$ for each $0 \leq i \leq p-1$. For each $p \in \mathbb{Z}^+$, the free vector space on the allowed $p$-paths on $(V,E)$ is denoted $\mathcal{A}_p=\mathcal{A}_p(G) = \mathcal{A}_p(V,E,\mathbb{K})$ and is called the \textit{space of allowed $p$-paths.} Furthermore, $\mathcal{A}_{-1}:=\mathbb{K}$, $\mathcal{A}_{-2}:= \{0\}$.

\textbf{Step 4:} The allowed paths do not form a chain complex since the image of an allowed path under $\partial$ need not to be allowed. To handle this problem, they define the \textit{space of $\partial$-invariant $p$-paths on $G$} as the following subspace of $\mathcal{A}_p(G)$:
$$
\Omega_p=\Omega_p(G)=\Omega_p(V,E,\mathbb{K}):=\{c \in \mathcal{A}_p: \partial_p(c) \in \mathcal{A}_{p-1}\}
$$
One further defines $\Omega_{-1}:=\mathcal{A}_{-1}=\mathbb{K}$ and $\Omega_{-2}:=\mathcal{A}_{-2}=\{0\}$. Now we have a chain complex and this yields to path homology groups.

\textbf{Filtration:} For any $\delta \in \mathbb{R}$ and a directed weighted graph $G=(V,E)$ with the edge weight function $W: V \times V \rightarrow \mathbb{R}$, we define the directed subgraphs $\mathfrak{G}_{\delta}=(V, E_{\delta})$ where $E_{\delta}:=\{(v,v') \in V \times V : v \neq v', W(v,v')\leq \delta \}$. We then define the \textit{digraph filtration} as 
$$\{G_{\delta} \hookrightarrow G_{\delta'} \}_{\delta \leq \delta'}.$$

After getting the filtration, instead of the persistent homology, they apply the persistent path homology (PPH). They show that in an undirected graph, PPH and Dowker filtrations agree in dimension 1 if a certain local condition is satisfied (need to be square-free). The authors also prove the stability result for PPH. 

\subsection{Generalizations of Vietoris-Rips filtration (GVR)} The author in \cite{turner2016generalizations} uses ordered-tuple complexes instead of using simplicial complexes to increase the flexibility regarding order. An ordered-tuple complex (OT-complex) is a collection $K$ of ordered tuples such that if $(v_0,v_1,v_2,...,v_n) \in K$ then $(v_0,...,\hat{v_i},...,v_n) \in K$ for all $i$ (where $(v_0,...,\hat{v_i},...,v_n)$ is the ordered tuple with $v_i$ removed). For example, the tuples $(v_1,v_2,v_3)$ and $(v_3,v_1,v_2)$ are distinct. One can define a chain complex, homology and $k$-th dimensional ordered tuple persistence homology of OT-complexes as defined for simplicial complexes. Then, the author defines the four generalization of Vietoris-Rips filtrations. She also proves the stability theorem for each case. Hence, we will not mention the theorem for each case separately. For the following filtrations, let $(V,f)$ be the vertex set and a function $f:V\times V \rightarrow \mathbb{R}$ ($f$ can be considered as an edge weight function).

\hypertarget{gvra}{\subsubsection{Vietoris-Rips filtration under sym$_a$}} For any $a \in [0, 1]$ we can define a symmetric function sym$_a(f):V\times V \rightarrow \mathbb{R}$ 
$$(u,v) \rightarrow a \text{min}\{f(u,v),f(v,u)\}+ (1-a)\text{max}\{f(u,v),f(v,u)\}.$$ 
Set $\mathcal{R}(V,\text{sym}_a(f))_t$ be the simplicial complex containing $[v_0,v_1,...,v_p]$ whenever sym$_af(v_i,v_j)\leq t$ for all $i,j$. This filtration is called the \textit{Vietoris-Rips filtration under sym$_a$} of $(V,f)$.

\hypertarget{gvrb}{\subsubsection{Directed Vietoris-Rips filtration}} Set $\{\mathcal{R}^{dir}(V,f)_t\}$ to be the filtration of OT-complexes where $(v_0,v_1,...,v_p)\in \mathcal{R}^{dir}(V)_t$ when $f(v_i,v_j)\leq t$ for all $i\leq j$. This filtration is called $\{\mathcal{R}^{dir}(X)_t\}$ the \textit{directed Vietoris-Rips filtration} of $(V,f)$.

\hypertarget{gvrc}{\subsubsection{Associated filtration of directed graphs}} There is a natural filtration of directed graphs $\{\mathcal{D}(V)_t: t \in [0,\infty )\}$ associated to $V$ by setting $\mathcal{D}(V,f)_t$ to the the directed graph with vertices $\{v \in V:f(v,v)\leq t \}$ and including the directed edges $u\rightarrow v$ whenever max$\{f(u,u), f(v,v), f(u,v)\} \leq t$. This filtration is called the \textit{associated filtration of directed graphs} of $(V,f)$.

 \hypertarget{gvrd}{\subsubsection{Preorder filtration}} Given a preorder $(V,\leq_t)$, let $\mathcal{O}(V,\leq_t)$ be the OT-complex containing $(v_0,v_1,...,v_p)$ when $v_0\leq v_1 \leq ... \leq v_p$. Let $\mathcal{O}(V,f)) =\{\mathcal{O}(V,f)_t\}$ be the filtration of OT-complexes corresponding to the filtration of posets $\{(V_t,\leq_t)\}$. This filtration is called the \textit{preorder filtration} of $(V,f)$

\begin{table}[t!]

\caption{The comparison of filtrations}
 \centering
\begin{tabular}{|c|c c c c c c c c c|}
\hline    & \rotatebox{90}{Undirected network} & \rotatebox{90}{Directed network} & \rotatebox{90}{Weighted network} & \rotatebox{90}{Unweighted network}& \rotatebox{90}{Metric network} & \rotatebox{90}{Dynamic network} &  \rotatebox{90}{Simplicial Complex} & \rotatebox{90}{Ordered-tuple complex} & \rotatebox{90}{Stable} \\ 

\hline    \hyperlink{vr}{VR} & \cmark  &  & \cmark &  &  &  &\cmark &  & \cmark  \\ 
\hline    \hyperlink{dss}{DSS} &  & \cmark & \cmark &  &  & & \cmark &  & \cmark \\
\hline    \hyperlink{ccl}{CCL} & \cmark &  &  & \cmark &  & & \cmark &  &  \\ 

\hline    \hyperlink{vbcl}{VBCL} & \cmark  & & \cmark &  &  & & \cmark &  &  \\ 
\hline    \hyperlink{kcl}{kCL} & \cmark &  & \cmark &  &  & & \cmark &  &  \\ 
\hline    \hyperlink{wscl}{WS} & \cmark &  & \cmark &  &  & & \cmark & &  \\ 
\hline    \hyperlink{vfb}{VFB} & \cmark &  & \cmark &  &  & & \cmark & & \cmark \\ 
\hline    \hyperlink{ic}{IC} & \cmark &  &  & \cmark & \cmark & & \cmark &  &  \\ 
\hline    \hyperlink{fmg}{FMG} & \cmark &  & \cmark &  & \cmark & & \cmark &  & \cmark \\
\hline    \hyperlink{pow}{POW} & \cmark &  &  & \cmark &  & & \cmark & &  \\
\hline    \hyperlink{tmp}{TMP} & & \cmark &  & \cmark &  & \cmark & \cmark & &  \\
\hline    \hyperlink{zsf}{ZSF} & \cmark &  &  & \cmark &  & \cmark & \cmark & & \cmark \\
\hline    \hyperlink{pph}{PPH} &  & \cmark & \cmark &  &  & & \cmark &  & \cmark \\ 
\hline    \hyperlink{gvra}{GVR-a} & \cmark &  & \cmark &  &  & & \cmark & & \cmark  \\
\hline    \hyperlink{gvrb}{GVR-b} &  & \cmark & \cmark &  &  &  & & \cmark & \cmark  \\
\hline    \hyperlink{gvrc}{GVR-c} &  & \cmark & \cmark &  &  & & \cmark &  & \cmark  \\
\hline    \hyperlink{gvrd}{GVR-d} &  & \cmark & \cmark &  &  &  & & \cmark & \cmark  \\
\hline
\end{tabular}
\label{table:conf}
\end{table}

\section{Algorithms and Applications}\label{sec:app}
In recent years, persistent homology has found applications in data analysis, including neuroscience \cite{sizemore2018importance}, time series data \cite{seversky2016time}, text mining \cite{wagner2012computational} and shape analysis \cite{gamble2010exploring}. In the complex network setting, while some studies analyze the evaluation of a single graph, some studies analyze multiple graphs for graph matching and classification with characterizing the temporal changes in topological features of a network. Besides these, while some studies use Betti numbers, some studies use persistent diagrams to extract some statistical features of the network. In this section, we categorize the persistent homology enabled applications as single graph and multiple graph analysis. We explain the algorithms and applications of each study in their corresponding sections. We also provide a comparison table, Table~\ref{table:appconfs} and Table~\ref{table:appconfm}, for algorithms with datasets after each section.

\subsection{Analysis on Single Graph}

In some applications, persistent homology is used to detect global structural features of a single network such as complexity and distributions of strongly connected regions. While some applications analyze the evaluation of a single graph according to edge weights, others analyze the evaluation of the graph over time. 

In many studies, Betti numbers are used as the complexity measure for different networks. Benzekry et al. \cite{benzekry2015design} propose that cancer therapy can be guided by changes in the complexity of protein-protein interaction (PPI) networks. They analyze 11 cancer interaction networks and find out that there is a correlation between 1-dimensional Betti number and survival of cancer patients. They compute Betti numbers using the power filtration (\hyperlink{pow}{POW}). To examine the effect of a node on the network complexity, each node in the network is removed and the change in Betti number is recorded. They consider the drop of the Betti number as the drop of the complexity. Therefore, if the removal of a node results in the largest drop in Betti number, it also results in the largest drop in complexity and is potentially a good drug target. 

Similar to this, Rucco et al.~\cite{rucco2016characterisation} use Betti number as persistent entropy to measure the graph complexity. They study the behavior of the idiotypic network of the mammal immune system. Their main goal is to detect the behavior of the immune system reaction to an external stimulus in terms of phase transitions. In addition to the persistent entropy, they use 2 other graph complexity measures, which are the connectivity entropy and the approximate von Neumann entropy~\cite{Petz2001}. While connectivity entropy is used to analyze the structural properties and to identify the set of key players of the idiotypic network, approximate von Neumann entropy is used to distinguish graphs corresponding to the same system but in different conditions. For persistent entropy, they use persistent barcodes constructed with the Vietoris-Rips filtration (\hyperlink{vr}{VR}). In their experiment, they create the simulation of the idiotypic network and obtain a weighted idiotypic network using the coexistence function as a weight function between antibodies. After computing the 3 different entropy measures on this network, they identify that peak on entropy corresponds to the activation of the immune response. While the connectivity entropy does not distinguish between the activation and the immune memory states, both the approximated von Neumann entropy and the persistent entropy are able to recognize the activation of the immune system. The analysis of the Betti numbers reveals that there is a subset of antibodies arranged in a 1-dimensional hole that is present both in the activation state and in the memory state.

Cliques and cycles are important structural features of complex networks to describe their cohesive structures. Rieck et al.~\cite{Rieck2018CliqueCP} use persistent homology to detect clique communities and their evolution in weighted networks. Persistent diagram is created using the vertex-based clique filtration (\hyperlink{vbcl}{VBCL}) and the $k$-clique filtration (\hyperlink{kcl}{kCL}). They analyze the connectivity relations for all clique degrees and all weight thresholds. Various networks are studied including co-occurrence network, brain network, and collaboration network. An interactive visualization tool is created that is capable of detecting and tracking the evolution of networks' clique communities for different thresholds and clique degrees.

Persistent homology is also used to analyze the brain networks by computing distributions of cliques (brain regions) and cycles (strongly connected regions) in them. In \cite{giusti2016two}, the authors review the underlying mathematical background of using simplicial complex in neural data, specifically brain networks. They list different types of simplicial complexes for encoding neural data such as networks, clique complex, independence complex, and concurrence complex. They also elaborate on using persistent homology to measure the global structure of simplicial complexes and the strength of neural connections using the weighted simplex filtration (\hyperlink{wscl}{WS}) to generate persistent diagrams.  

In \cite{sizemore2018cliques}, the authors test the hypothesis that the spatial distributions of cliques and cycles will differ in their anatomical locations. They construct 1- and 2-dimensional persistent diagrams of brain networks using the Vietoris-Rips filtration (\hyperlink{vr}{VR}).  The structural brain networks of eight volunteers is extracted using diffusion spectrum imaging. The undirected and weighted network consists of 83 nodes representing different brain regions and edges that refer to the density of white matter between the nodes. Weak and strong connections between cliques are assessed by observing the difference between birth and death times of $k$-cliques in persistent diagrams. 

Additionally, in \cite{chung2013persistent}, persistent homology is also used to analyze the brain networks with the aim of examining the abnormal white matter in maltreated children. Networks are obtained by thresholding (based on the sample covariance) sparse correlations for the Jacobian determinant from magnetic resonance imaging (MRI) and fractional anisotropy from diffusion tensor imaging (DTI) at different threshold values. The collection of the thresholded graphs forms a Vietoris-Rips filtration (\hyperlink{vr}{VR}).

Moreover, in \cite{khalid2014tracing}, the authors demonstrate that persistent homology is useful in analyzing functional brain connectivity. The application involves electroencephalography (EEG) data from eight cortical regions of corticosterone (CORT) induced depression mouse and control models. After the EEG measurement is obtained, the square root of (1-correlation) distance metric is used to create a binary network. Next, the Vietoris-Rips filtration (\hyperlink{vr}{VR}) is applied and used to visualize topological changes by 0-dimensional barcodes which are then used to construct single-linkage dendrograms (SLD). Finally, single-linkage distance is computed using the generated SLDs. The results show that CORT model is characterized by an increased local connectivity and by a decreased global connectivity.

Besides its utility on brain networks, persistent homology is also used to analyze word co-occurrence, remittance, and migration networks. In \cite{salnikov2018co}, the authors study the word co-occurrence networks to explore the conceptual landscape of mathematical research. They first create the network using 54177 articles in arXiv from 01/1994 to 03/2007. Then they parse a concept list from Wikipedia that includes 1612 equations, theorems, and lemmas. Next, they combine these two datasets by checking 1612 concepts' appearance in the articles and find that 1067 of them match in at least one article and 35018 articles contain at least one of the concepts. They first take 1067 concepts as nodes and include a $(n-1)$-simplex for each article containing $n$-concepts. Furthermore, whenever the concept sets of two articles intersect at $n$ concepts, their corresponding simplices share a face of dimension $(n-1)$. In total, this construction results in 32707 unweighted edges. They use the temporal filtration (\hyperlink{tmp}{TMP}) using article dates. They create the 1- and 2-dimensional persistent diagrams, i.e. they just look at the 2-dimensional holes bounded by edges and 3-dimensional holes bounded by triangles respectively. They interpret these holes to explain the intrinsic characteristics of how research evolves in mathematics. They also explore the authors' conceptual profile using the holes and their attributes to the holes.

Ignacio et al. \cite{Ignacio2019} analyze the patterns and shapes in remittance and migration networks as a directed weighted network via persistent homology to identify flow patterns between multiple countries. They detect both local and global patterns that highlight simultaneous interactions among multiple nodes. They extend the Vietoris-Rips filtration (\hyperlink{vr}{VR}) to detect topological features such as persistent cycles in directed networks using the weight of the edges and create persistence barcodes. They use 0-, 1- and 2-dimensional barcodes to analyze the cycles in networks. As a modification on 1 and 2-dimensional barcodes, to encode additional information, they color the bars in barcodes according to the standard deviation of the weights in the cycles they represent. They create the 2015 Asian net migration and remittance networks which include 50 countries and states to perform their analysis on. They define the weight of a directed edge $(a,b)$ as the profit country $b$ gains from exchanging remittances with the country $a$ for remittance networks and define it in a similar manner for net migration networks.

One of the challenges for most graphing methods is the inability of visualizing the global structure of graphs as a result of the absence of interactive exploration mechanisms. Persistent homology is used to address this challenge~\cite{DBLP:journals/corr/abs-1712-05548}. They use 0-dimensional PH features to control and modify force-directed layouts of a graph. The 0-dimensional barcode, obtained by the power filtration (\hyperlink{pow}{POW}), enables the visualization of contraction and repulsion events in the network. More forces are added to the graph layout based on the selected number of barcodes. They have three case studies to show the effectiveness of their method on 3 different real-world networks.  One of the networks is ``Les Miserables" which contains 77 nodes (characters) and 254 edges, weighted by how many scenes two characters share during any chapter of the novel. Some of the key characters featured in the book can be identified on the force-directed layout modified with PH features. They are also able to extract major important nodes in the Madrid Train Bombing network and US Senate 2007 and 2008 Co- and Anti-voting network using their method.

\begin{table}[t!]

\caption{The comparison of algorithms and applications on Single Graph}
 \centering
\resizebox{\textwidth}{!}{\begin{tabular}{|c| c| c| c | c |c |} 
\hline    Paper & Filtration & Topological Summary  & Data \\
\hline
\hline\cite{benzekry2015design} & \hyperlink{pow}{POW}  & 0-1 dim Betti numbers & PPI networks\\ 

\hline\cite{sizemore2018cliques} & \hyperlink{vr}{VR}  & 1-2 dim PD  & Brain networks\\
\hline \cite{chung2013persistent} & \hyperlink{vr}{VR}  & 0 dim PD  & Brain networks\\

\hline \cite{giusti2016two} & \hyperlink{wscl}{WS}  & 0-2 dim PD & Brain Networks\\

\hline \cite{khalid2014tracing} & \hyperlink{vr}{VR}  & 0-dim PD & Brain Networks\\

\hline \cite{Ignacio2019} &  \hyperlink{vr}{VR} & 0-2 dim PD & Migration and remittance networks  \\
\hline\cite{rucco2016characterisation} & \hyperlink{vr}{VR}  & 1 dim PB  & Simulated idiotypic networks\\

\hline \cite{salnikov2018co} &  \hyperlink{tmp}{TMP} & 1-2 dim PD & Co-occurrence networks\\
\hline\cite{DBLP:journals/corr/abs-1712-05548} & \hyperlink{pow}{POW}  & 0 dim PB & {Co-occurrence networks} \\
\hline \cite{Rieck2018CliqueCP} & \hyperlink{vbcl}{VBCL},\hyperlink{kcl}{kCL}   & 0 dim PD  & Co-occurrence, brain and collaboration network \\

\hline 
\end{tabular}}
\label{table:appconfs}
\end{table}

\subsection{Multiple Graphs Analysis}
Graph comparison is an important task for many graph applications such as classification and matching. On the other hand, it is a computationally complex problem where we need to compute the similarity between 2 networks~\cite{conte2004thirty}. It has been studied for many years and defined as either exact matches (e.g. graph isomorphism ~\cite{cordella2004sub}) or some measures of structural similarity (e.g. graph edit distance~\cite{gao2010survey}). Graph kernels are also used to capture the graph similarity~\cite{shervashidze2011weisfeiler}. Recent years, persistent homology is used to extract topological features of networks to compare them. 

While most existing metrics for network structure rely on local features of vertices such as node degrees, correlation of neighborhood nodes, they do not capture the precise mesoscopic structure of complex networks. Sizemore et al.~\cite{sizemore_classification} extract mesoscale homological features as 0-3 dimensional Betti numbers. They use the Vietoris-Rips (\hyperlink{vr}{VR}) filtration to compute the homology and record the maximal clique distribution and Betti sequence. Extracted features are used to classify 14 commonly studied weighted network models into four groups or classes with agglomerative hierarchical clustering to use for graph classification. Betti values and parameters from the maximal clique distribution are used to determine the structural similarities between networks. After classifying networks into groups, they analyze the structural patterns in each group of networks.

In \cite{petri2013topological, petri2013networks, binchi2014jholes}, persistent homology is used to detect particular non-local structural features of networks. After creating the barcodes with the inverse Vietoris-Rips (\hyperlink{vr}{VR}) filtration based on edge weights, statistical distributions of 1-dimensional barcodes are computed. They classify real-world networks into 2 classes according to the similarity of their cycle distribution with randomized version. In Class I, cycle distributions are markedly different from the randomized versions and in Class II, cycle distributions are very close to their random versions. The authors study different network datasets, such as US air passenger networks, C. Elegans's neuronal network \cite{watts1998collective}, the online messages network \cite{opsahl2009clustering}, gene network, network of mentions and re-tweet between Twitter users, school face-to-face contact network, co-authorship networks. While the gene network and airport network are in class 1, co-authorship networks and twitter network are in class 2. 

In \cite{carriere2019general}, the authors propose to use persistence diagrams for graph classification problem for undirected weighted graphs. They first define a graph kernel function, namely heat kernel signatures \cite{hu2014stable}, on networks and use the sublevel and superlevel \hyperlink{vfb}{VFB} filtration on each network to generate PDs. Then they employ two layers neural network architecture to process the PDs and classify the graphs. They evaluate their classification model on social networks, medical and biological networks. They also compare their results with four different state-of-the-art graph classification methods and show that their method has comparable results despite being much simpler than other methods. 

Moreover, persistent homology is used to analyze the structure of weighted networks. In \cite{carstens_horadam_2013}, the authors consider the collaboration networks as weighted network. They use the Vietoris-Rips (\hyperlink{vr}{VR}) filtration to generate the persistence barcodes of networks. They employ the Betti numbers of 0, 1, and 2 dimensions and use them to distinguish collaboration networks from random networks. They conclude that the first and second Betti numbers give us richer information about weighted networks.
 
Siddharth et al.~\cite{Siddharth2017} study the growing collaboration network with a temporal parametrization and characterize the temporal changes in its topological features. In a collaboration network, each person in a paper or a movie is represented as a vertex, and each collaborative act (and each of its subsets) is represented as a simplex of vertices comprising it. They define a temporal filtration (\hyperlink{tmp}{TMP}) from growing collaboration networks, with adding new collaborations occurred in each year. In addition, they introduce a new distance measure between a growing network which captures the difference in the rate of growth of cycles in the networks being compared. They use DBLP (Digital bibliography \& library project) and IMDB (Internet movie database) data sets from 1950-2008 considering 10-year windows. They study the topological properties of networks as the growth in the cyclicity, with respect to the time corresponding to the 10-year windows, and size of the largest connected component. 

In \cite{Schauf2016DiscriminationOE}, the authors consider the national input-output networks of domestic products as a weighted network and use persistent homology to identify dissimilarities between them. The nodes are available sectors in an economy and edge weights are the monetary flow measuring the magnitude of the economic relationship between two sectors. They generate persistence diagrams for dimensions 0, 1, and 2 with the Vietoris-Rips (\hyperlink{vr}{VR}) filtration. Using 0-dimensional diagrams, they distinguish economies with high GDP, large population, and small import/export percentages of GDP from those with lower GDP, small population, and larger import/export percentages. They also discuss the potential for applying higher-dimensional persistent homology to study these networks.

Similarly, financial networks are considered as weighted networks and persistent homology is used to
detect early signs of critical transitions of financial crisis in~\cite{gidea2017topological}. The vertices correspond to the stocks, each pair of distinct nodes is connected by an edge and each edge is assigned a weight using the Pearson correlation coefficient. For each time frame, they generate 0- and 1-dimensional persistent diagrams of the network using the Vietoris-Rips (\hyperlink{vr}{VR}) filtration. Then, the distance between them is measured via Wasserstein distance. They show that the persistent diagrams and the distances between them have significant changes prior to the 2007-2008 financial crisis. 

Furthermore, in \cite{keil2018topological}, the undirected attributed networks are considered as weighted networks. They first assign weights on edges using the vertex attributes. Then, they extract the ego-networks of each vertex and define a graph kernel function, namely the diffusion Fr{\'e}chet function \cite{martinez2018probing}, on each ego network that takes both the network topology and edge weights into consideration. Next, they generate the persistence diagrams of each ego network using the sublevel and the superlevel \hyperlink{vfb}{VFB} filtration and obtain the distance matrix between each vertex computing the Wasserstein distance between their persistence diagrams. Finally, they cluster the network using the $k$-means clustering algorithm.  

Beside previous weighted networks, brain networks are considered as sparse weighted networks and persistent homology is also used to analyze them~\cite{chung2015persistent}. They obtain the topological structure of a graph induced by sparse correlation. They first transform MRI and DTI data to weighted networks where they employ the sparse Pearson correlation to obtain the edge weights. They generate the 0-dimensional Betti plots for the brain networks using the Vietoris-Rips filtration (\hyperlink{vr}{VR}). They also generate Betti plots using sparse covariance. They show that the sparse correlation method gets a huge group separation between normal and stress-exposed children visually. This method is also less computationally expensive than the sparse covariance method. 

In \cite{knyazeva2018resting}, the authors study dynamical connectome state analysis on brain networks using three different methods: $k$-means clustering, modularity based clustering and topological feature based clustering. They consider brain networks as weighted networks. In topological feature based clustering, they use the Vietoris-Rips (\hyperlink{vr}{VR}) filtration. They first split the correlation matrix to the two matrices with positive and negative correlations. Then, they create \hyperlink{vr}{VR} filtrations for both matrices. In their clustering, different type of connections describes different processes in the brain, so they compute persistent homology with annotated intervals collection. After getting the intervals, they compute different statistics for each homology group and for types of interactions. Then, finally, they perform hierarchical clustering based on these topological features. They show that topological feature based clustering is more informative than the other two clustering methods.

In addition to this, in \cite{chowdhury2016persistent,chowdhury2018persistent}, the authors classify the brain (hippocampal) networks using persistence diagrams. They consider five different environments with $4$ holes, $3$ holes, $2$ holes, $1$ hole and no hole and for each environment, 20 simulated brain networks are created. Persistence diagrams of these 100 networks are computed with the Dowker filtration (\hyperlink{dss}{DSS}). They use the bottleneck distance between the 1-dimensional diagrams of networks to compare them. Finally, they classify the networks using the single linkage dendrogram algorithm and show that Dowker filtration is successful in capturing the differences between the five classes of networks. The authors also work on the same problem and dataset using the zigzag simplicial filtration (\hyperlink{zsf}{ZSF}) in \cite{chowdhury2018importance}. They create 1-dimensional zigzag persistent diagrams to perform persistent homology computations on dynamic simplicial complexes resulted from these brain networks.

In \cite{yoo2016topological}, the authors show that persistent homology, or more precisely persistence vineyard, is a robust approach to estimate functional connectivity in the resting and gaming stages of the brain networks. They conduct an experiment with 26 male college students aged 19-29 years old from two universities located in Seoul, Republic of Korea. They undergo all the 26 healthy subjects resting and gaming experiments. Each stage was recorded for five minutes separately. They segment their data using 30s window lengths and 2s step size. For each window, they compute the persistence diagram using Pearson correlation between brain channels employing the weighted simplex filtration (\hyperlink{wscl}{WS}). Then, they compute the 0-dimensional persistence vineyard to analyze the dynamic brain connectivity. In a brief, a persistent vineyard is a $p$ dimension persistent diagram with a time dimension added, tracking the birth and death of $p$ dimension diagrams in a time-varying topological space \cite{edelsbrunner2010computational}. Their results show that persistence vineyard is successful to determine the temporarily dynamic properties of the brain in a robust and threshold-free way. They also show that persistent vineyard is more effective than the principal component analysis (PCA) and standard graph theoretical methods.

\cite{petri2014homological} compares resting state functional brain activity in 15 healthy volunteers after intravenous infusion of placebo and psilocybin using persistent homology and other statistical methods-density function. First, the raw data from fMRI (functional magnetic resonance imaging) dataset is transformed into a functional network. They create the 1-dimensional persistence diagram using the Vietoris-Rips (\hyperlink{vr}{VR}) filtration. Later, they define two different homological scaffolds depending on how frequently edges are part of the generators of the persistent homology groups and how persistent are the generators to which they belong to. The results show that the homological structure of the brain's functional patterns undergoes a dramatic change post-psilocybin, characterized by the appearance of many transient structures of low stability and of a small number of persistent ones that are not observed in the case of placebo.

In \cite{horak2009persistent}, the authors first study random graphs using the clique (\hyperlink{ccl}{CCL}) filtration. Using different probabilities, they generate random networks and compute their barcodes. They show that the results on these barcodes are in agreement with the theoretical studies on these complexes. As another application, they study an email network. They create barcodes and show that higher dimensional barcodes, which do not exist for random networks, correspond to more dense communications among certain groups. They also apply their methods on scale-free networks with a modular structure. They use three different parameters to generate three types of scale networks: Clustered modular networks, clustered non-modular networks and non-clustered modular networks. They show that both clustered modular and clustered non-modular networks have more bars in their 3-dimensional and 4-dimensional barcodes than the non-clustered modular network. 

Moreover, persistent homology is used for metric graph comparison. In \cite{dey2015comparing}, the authors first introduce the functional metric graph filtration (\hyperlink{fmg}{FMG}) on metric graphs. Then, they define the persistence distortion distance between two finite metric graphs using the persistence diagrams from FMG filtration. In their experiment,  they show the stability of the proposed distance measure on the Athen's road network as a metric graph and generate its noisy sample using a noise level $\epsilon$. The results show that the persistence distortion distance between the original graph and its noisy sample grows roughly proportionally to $\epsilon$. They also use proposed persistence distortion distance to compare surface meshes of different geometric models. Models from the same group have very smaller persistence distortion distances among them than those between the dissimilar group, which shows the that proposed distance is able to differentiate surface models.

In \cite{Hajij2018VisualDO}, persistent homology is employed to quantify structural changes in time-varying (dynamic) graphs. Their objective is to transform each instance of the time-varying graph into a metric space, extract topological features using persistent homology, and compare those features over time by means of bottleneck or Wasserstein distance between their corresponding persistence diagrams. Finally, several case studies on real-world networks, such as high school communication network, show how this method can find cyclic patterns, deviations from those patterns, and one-time events in time-varying graphs. In particular, 0- and 1-dimensional PH are utilized to detect the components and tunnels respectively. Each graph constitutes a distinct metric space for which the Vietoris-Rips (\hyperlink{vr}{VR}) filtration is implemented to compute the corresponding Betti numbers.

High order networks are weighted complete hypergraphs collecting relationships between elements of tuples. Computing distance between high order network is difficult  when the number of nodes is large. In \cite{huang2017persistent}, the authors use persistent homology to derive distance approximations of networks. They compute the bottleneck distance between persistence diagrams of networks to evaluate the differences between networks. They first define a relationship function between a set of nodes to represent a measure of similarity or dissimilarity for members of the group. They use this function to assign weight on each simplex. Using these weights and the weighted simplex filtration (\hyperlink{wscl}{WS}), they generate 0-, 1- and 2-dimensional persistence diagrams.  They show that they can lower bound distance between two higher order networks, which is in general computationally expensive, with a computationally less expensive distance between their persistence diagrams. They apply their method to the coauthorship networks. They first create the networks using 5 journals from the mathematics community and 6 journals from the engineering community. They use the lower bounds to classify the networks, distinguish the collaboration patterns of engineering and mathematics community and also discriminate engineering communities with different research interests.

To answer the question of whether the existing anonymization mechanisms for preserving privacy truly keep the graph utility, Gao et al.~\cite{gao2018studying} employ persistent homology to analyze and evaluate four anonymization mechanisms. They study online social networks (OSN). They define the distance between two nodes as the number of hopes on the shortest path between these nodes and create 0-,1- and 2-dimensional persistence barcodes using this distance in the power filtration (\hyperlink{pow}{POW}). They analyze the original and anonymized OSNs using the barcodes. The results show that original OSN graphs have stable structures. Furthermore, the 0-dimensional barcodes they obtain show that most anonymized OSNs are more closely connected than the original graph. All anonymized graphs are not as stable as the original graph, because they have more 2-dimensional holes or larger holes. They also compare their results with traditional graph metrics. 

In \cite{kim2017stable}, the authors study flocking/swarming behaviors in animals. They first create dynamic graphs and simplicial complexes using the Vietoris-Rips complex for a fixed scale parameter. Then, they construct the zigzag simplicial filtration (\hyperlink{zsf}{ZSF}) and obtain 0-dimensional zigzag persistent diagrams to classify the four different type of flocking behavior in animals. Finally, using the bottleneck distance, the single linkage hierarchical clustering, and MDS, they distinguish the 4 behaviors very well.

In \cite{chowdhury2018persistent2}, the authors characterize the directed cycle networks by digraph filtration using persistent path homology (\hyperlink{pph}{PPH}). They prove that the persistent diagrams of a cycle network with $n$ nodes for $n\geq 3$ solely depends on $n$. 

\begin{table}[t!]

\caption{The comparison of algorithms and applications on Multiple Graphs}
\centering
\resizebox{\textwidth}{!}{\begin{tabular}{|c| c| c| c | c |c |} 
\hline    Paper & Filtration & Topological Summary & Data \\
 \hline

\hline\cite{chowdhury2016persistent, chowdhury2018persistent} &   \hyperlink{dss}{DSS} & 1 dim PD & Brain Networks\\

\hline \cite{chowdhury2018importance} &  \hyperlink{zsf}{ZSF}  & 1 dim zigzag PD & Brain Networks\\

\hline \cite{chung2015persistent} & \hyperlink{vr}{VR}  & 0 dim Betti plots & Brain networks  \\
\hline\cite{knyazeva2018resting} &  \hyperlink{vr}{VR} & 0-1 dim PD & Brain networks\\
\hline \cite{petri2014homological} & \hyperlink{vr}{VR}  & 1 dim PD & Brain networks  \\
\hline \cite{yoo2016topological} & \hyperlink{wscl}{WS}  & 0 dim persistence vineyards & Brain Networks\\

%\hline \cite{costamethodology} & \hyperlink{vr}{VR}  & 0- and 1-dim PD & Sparse brain networks  \\

\hline\cite{huang2017persistent} &  \hyperlink{wscl}{WS} & 0-2 dim PD & Collaboration networks\\

\hline\cite{carstens_horadam_2013} &  \hyperlink{vr}{VR} & 0-2 dim Betti numbers & Collaboration networks\\

\hline \cite{Siddharth2017} & \hyperlink{tmp}{TMP}  & 1 dim Betti numbers& Collaboration networks  \\

\hline \cite{sizemore_classification} & \hyperlink{vr}{VR}  & 0-3 dim Betti numbers &  PPI, brain and simulated weighted networks \\
\hline \cite{Schauf2016DiscriminationOE} & \hyperlink{vr}{VR}  & 0-2 dim PD &  Economy networks\\

\hline \cite{gidea2017topological} & \hyperlink{vr}{VR}  & 0-1 dim PD & Finance networks\\

\hline\cite{horak2009persistent} &  \hyperlink{cl}{CL} & 0-11 dim PD & Random, email and scale-free networks\\

\hline\cite{dey2015comparing} &  \hyperlink{fmg}{FMG} & 0 dim PD & Road networks\\
\hline\cite{kim2017stable} &  \hyperlink{vsf}{VSF} & 0 dim zigzag PD & Dynamic biological networks\\
\hline\cite{chowdhury2018persistent2} &  \hyperlink{pph}{PPH} & 1 dim PD & Cycle networks\\
\hline \cite{Hajij2018VisualDO} & \hyperlink{POW}{POW}  & 0-1 dim PD &  Dynamic communication networks\\
\hline \cite{keil2018topological} & \hyperlink{vfb}{VFB}  & 0-1 dim PD & Attributed social networks\\

\hline \cite{gao2018studying} & \hyperlink{pow}{POW}  & 0-1 dim PD & Online social networks  \\
\hline \cite{carriere2019general} & \hyperlink{vfb}{VFB}  & 0-1 dim PD & Social, medical and biological networks\\
\hline \cite{petri2013topological, petri2013networks, binchi2014jholes} &  \hyperlink{vr}{VR}& 1 dim PB & Social, infrastructural and biological networks\\
\hline 
\end{tabular}}

\label{table:appconfm}
\end{table}

\section{Conclusion}\label{sec:conc}

In this paper, we provide a conceptual review of key advancements in the area of using PH on applied network science. We look into research studies that use PH on networks and highlight different algorithms that are used to extract topological features of networks.  We review the applications where PH is used in solving network mining problems. We believe our summary of the analysis of PH on networks will provide important insights to researchers in applied network science.

At the moment, the implicit goal of most studies is to extract the topological features of the networks that persist across multiple scales. However, there are some limitations to these studies. Firstly, scalability may be a concern for future progress. The networks used in these studies are mostly small networks (number of vertices is less than 1000). There is still significant work needed to be done in scaling PH approaches for larger networks. Secondly, there are some filtrations whose stability has not proven yet. 

Furthermore, although many filtration methods are proposed, they are mainly designed for static networks. However, many real-world networks are evolving over time. For example, in the Facebook network, friendships between users always dynamically change over time, hence new edges are continuously added to the social network while some edges may be deleted. Most of the existing methods cannot be directly applied to large scale evolving networks. New filtration algorithms, which are able to tackle the dynamic nature of evolving networks, are highly desirable in persistent homology. 

As another future research direction, PH can also be used for network and sub-network embedding problems.      

\begin{backmatter}
\section*{Abbreviations}
PH, Persistent Homology; PB, Persistence Barcode; PD, Persistence Diagram; VR, Vietoris-Rips Filtration; DSS, Dowker Sink and Source Filtration; CCL, Clique Complex Filtration;  WRCL, Weight Rank Clique Filtration; VBCL, Vertex-Based Clique Filtration; kCL, $k$-clique filtration; WS, Weighted Simplex Filtration; VFB, Vertex function based filtration; IC, Intrinsic \v{C}ech Filtration; FMG, Functional Metric Graph Filtration; POW, Power Filtration; TMP, Temporal Filtration; ZSF, Zigzag Simplicial Filtration; PPH, Persistent Path Homology; GVR, Generalizations of Vietoris-Rips Filtration; PPI, Protein-protein interaction networks; MRI, Magnetic Resonance Imaging; fMRI, Functional Magnetic Resonance Imaging; DTI, Diffusion Tensor Imaging; EEG, Electroencephalography; CORT, Cortical Regions of Corticosterone; SLD, Single-linkage Dendrograms; US, United States; GDP, Gross Domestic Product; OSN, Online Social Networks; DBLP, Digital Bibliography \& Library Project; IMDB, Internet Movie Database; PCA, principal component analysis; MDS, multidimensional scaling; dim, dimension.

\section*{Availability of data and material}
Not applicable

\section*{Competing interests}
  The authors declare that they have no competing interests.

\section*{Funding}
Not applicable

\section*{Author's contributions}
MEA and EA designed the project. MEA studied the filtrations defined on networks. MEA, EA and AEF studied the algorithm and applications of the persistent homology in network settings. MEA, EA and AEF wrote the manuscript.
 %   Text for this section \ldots

\section*{Acknowledgements}
Not applicable
%  Text for this section \ldots
%%  After submission of the .TEX file,                     %%
%%  you will be prompted to submit your .BBL file.         %%

% if your bibliography is in bibtex format, use those commands:
\bibliographystyle{bmc-mathphys} % Style BST file (bmc-mathphys, vancouver, spbasic).
\bibliography{bmc_article}      % Bibliography file (usually '*.bib' )

\end{backmatter}
\end{document}